\documentclass[
	preprint
	number,
	sort&compress
]{elsarticle}

\usepackage[english]{babel}
\usepackage[utf8]{inputenc}
\usepackage{amsmath}
\usepackage{siunitx}
\usepackage{graphicx}
\usepackage{amssymb}
\usepackage{cases}
\usepackage{dsfont}
\usepackage{import}
\usepackage{mathtools}
\usepackage{xparse}
\usepackage{mleftright}
\usepackage[low-sup]{subdepth} %
\usepackage{textcomp} %
\usepackage{microtype}
\usepackage{changepage}
\usepackage{booktabs} %
\usepackage{bm}
\usepackage{hyperref}

\hypersetup{colorlinks=true,linktoc=all,linkcolor=blue,breaklinks=true,citecolor=blue,urlcolor=blue}
\allowdisplaybreaks[1]

\usepackage{flafter} %
\usepackage{wrapfig}

\makeatletter
\newcommand\RedeclareMathOperator{%
  \@ifstar{\def\rmo@s{m}\rmo@redeclare}{\def\rmo@s{o}\rmo@redeclare}%
}
\newcommand\rmo@redeclare[2]{%
  \begingroup \escapechar\m@ne\xdef\@gtempa{{\string#1}}\endgroup
  \expandafter\@ifundefined\@gtempa
     {\@latex@error{\noexpand#1undefined}\@ehc}%
     \relax
  \expandafter\rmo@declmathop\rmo@s{#1}{#2}}
\newcommand\rmo@declmathop[3]{%
  \DeclareRobustCommand{#2}{\qopname\newmcodes@#1{#3}}%
}
\@onlypreamble\RedeclareMathOperator
\makeatother

\usepackage{multirow}
\usepackage{xcolor}

\definecolor{uablue}{RGB}{0,68,102}
\colorlet{uablue100}{uablue}
\colorlet{uablue75} {uablue!75!white}
\colorlet{uablue50} {uablue!50!white}
\colorlet{uablue25} {uablue!25!white}
\colorlet{uablue10} {uablue!10!white}
\colorlet{uablue5}  {uablue!5!white}

\definecolor{uared}{RGB}{136,17,51}
\colorlet{uared100}{uared}
\colorlet{uared75} {uared!75!white}
\colorlet{uared50} {uared!50!white}
\colorlet{uared25} {uared!25!white}
\colorlet{uared10} {uared!10!white}
\colorlet{uared5}  {uared!5!white}

\definecolor{ualightblue}{RGB}{51,153,204}
\colorlet{ualightblue100}{ualightblue}
\colorlet{ualightblue75} {ualightblue!75!white}
\colorlet{ualightblue50} {ualightblue!50!white}
\colorlet{ualightblue25} {ualightblue!25!white}
\colorlet{ualightblue10} {ualightblue!10!white}
\colorlet{ualightblue5}  {ualightblue!5!white}

\definecolor{uagold}{RGB}{221,153,17}
\colorlet{uagold100}{uagold}
\colorlet{uagold75} {uagold!75!white}
\colorlet{uagold50} {uagold!50!white}
\colorlet{uagold25} {uagold!25!white}
\colorlet{uagold10} {uagold!10!white}
\colorlet{uagold5}  {uagold!5!white}

\definecolor{uayellow}{RGB}{170,170,0}
\colorlet{uayellow100}{uayellow}
\colorlet{uayellow75} {uayellow!75!white}
\colorlet{uayellow50} {uayellow!50!white}
\colorlet{uayellow25} {uayellow!25!white}
\colorlet{uayellow10} {uayellow!10!white}
\colorlet{uayellow5}  {uayellow!5!white}

\definecolor{darkgrey}{RGB}{60,60,59}
\definecolor{lightgrey}{RGB}{146,148,151}

\colorlet{dcolor}{uablue} %
\colorlet{dcolor100}{dcolor}
\colorlet{dcolor75} {dcolor!75!white}
\colorlet{dcolor50} {dcolor!50!white}
\colorlet{dcolor25} {dcolor!25!white}
\colorlet{dcolor10} {dcolor!10!white}
\colorlet{dcolor5}  {dcolor!5!white}

\colorlet{lcolor}{uared} %
\colorlet{lcolor100}{lcolor}
\colorlet{lcolor75} {lcolor!75!white}
\colorlet{lcolor50} {lcolor!50!white}
\colorlet{lcolor25} {lcolor!25!white}
\colorlet{lcolor10} {lcolor!10!white}
\colorlet{lcolor5}  {lcolor!5!white}

\colorlet{emphcolor}{uared}
\colorlet{emphcolor100}{emphcolor}
\colorlet{emphcolor75} {emphcolor!75!white}
\colorlet{emphcolor50} {emphcolor!50!white}
\colorlet{emphcolor25} {emphcolor!25!white}
\colorlet{emphcolor10} {emphcolor!10!white}
\colorlet{emphcolor5}  {emphcolor!5!white}

\colorlet{addcolor}{uagold} %
\colorlet{addcolor100}{addcolor}
\colorlet{addcolor75} {addcolor!75!white}
\colorlet{addcolor50} {addcolor!50!white}
\colorlet{addcolor25} {addcolor!25!white}
\colorlet{addcolor10} {addcolor!10!white}
\colorlet{addcolor5}  {addcolor!5!white}

\makeatletter
\DeclareFontFamily{OMX}{MnSymbolE}{}
\DeclareSymbolFont{MnLargeSymbols}{OMX}{MnSymbolE}{m}{n}
\SetSymbolFont{MnLargeSymbols}{bold}{OMX}{MnSymbolE}{b}{n}
\DeclareFontShape{OMX}{MnSymbolE}{m}{n}{
	<-6>  MnSymbolE5
	<6-7>  MnSymbolE6
	<7-8>  MnSymbolE7
	<8-9>  MnSymbolE8
	<9-10> MnSymbolE9
	<10-12> MnSymbolE10
	<12->   MnSymbolE12
}{}
\DeclareFontShape{OMX}{MnSymbolE}{b}{n}{
	<-6>  MnSymbolE-Bold5
	<6-7>  MnSymbolE-Bold6
	<7-8>  MnSymbolE-Bold7
	<8-9>  MnSymbolE-Bold8
	<9-10> MnSymbolE-Bold9
	<10-12> MnSymbolE-Bold10
	<12->   MnSymbolE-Bold12
}{}

\let\llangle\@undefined
\let\rrangle\@undefined
\DeclareMathDelimiter{\llangle}{\mathopen}%
{MnLargeSymbols}{'164}{MnLargeSymbols}{'164}
\DeclareMathDelimiter{\rrangle}{\mathclose}%
{MnLargeSymbols}{'171}{MnLargeSymbols}{'171}
\makeatother
	\DeclarePairedDelimiter{\opbrack}{(}{)}
	\NewDocumentCommand{\oplbrack}{s O{} m}{%
		\IfBooleanTF{#1}
		{\mleft(#3\mright.\kern-\nulldelimiterspace}
		{\mathopen{#2(} #3 \mathclose{}}%
	}
	\NewDocumentCommand{\oprbrack}{s O{} m}{%
		\IfBooleanTF{#1}
		{\kern-\nulldelimiterspace\mleft.#3\mright)}
		{\mathopen{} #3 \mathclose{#2)}}%
	}
	
	\NewDocumentCommand{\opsqlbrack}{s O{} m}{%
		\IfBooleanTF{#1}
		{\mleft[#3\mright.\kern-\nulldelimiterspace}
		{\mathopen{#2[} #3 \mathclose{}}%
	}
	\NewDocumentCommand{\opsqrbrack}{s O{} m}{%
		\IfBooleanTF{#1}
		{\kern-\nulldelimiterspace\mleft.#3\mright]}
		{\mathopen{} #3 \mathclose{#2]}}%
	}

	\NewDocumentCommand{\bracketlr}{s O{} m m m m m}{%
		\IfBooleanTF{#1}
		{\left#4 #3 \right#5^{#6}_{#7}}
		{\mathinner{ \mathopen{#2#4} #3 \mathclose{#2#5}^{#6}_{#7} } }%
	}
	\NewDocumentCommand{\bracketl}{s O{} m m}{%
		\IfBooleanTF{#1}
		{ \left#4 #3 \right.\kern-\nulldelimiterspace }
		{ \mathinner{ \mathopen{#2#4} #3 \mathclose{} } }%
	}
	\NewDocumentCommand{\bracketr}{s O{} m m m m}{%
		\IfBooleanTF{#1}
		{ \kern-\nulldelimiterspace\left. #3 \right#4^{#5}_{#6} }
		{ \mathinner{ \mathopen{} #3 \mathclose{#2#4}^{#5}_{#6} } }%
	}

	\NewDocumentCommand{\nbrack}{s O{} m t^G{} t_G{} }{%
		\IfBooleanTF{#1}
		{ \bracketlr*{#3}{(}{)}{#5}{#7} }
		{ \bracketlr[#2]{#3}{(}{)}{#5}{#7} }
	}
	\NewDocumentCommand{\nlbrack}{s O{} m}{%
		\IfBooleanTF{#1}
		{ \bracketl*{#3}{(} }
		{ \bracketl[#2]{#3}{(} }%
	}
	\NewDocumentCommand{\nrbrack}{s O{} m t^G{} t_G{} }{%
		\IfBooleanTF{#1}
		{ \bracketr*{#3}{)}{#5}{#7} }
		{ \bracketr[#2]{#3}{)}{#5}{#7} }%
	}

	\NewDocumentCommand{\sqbrack}{s O{} m t^G{} t_G{} }{%
		\IfBooleanTF{#1}
		{ \bracketlr*{#3}{[}{]}{#5}{#7} }
		{ \bracketlr[#2]{#3}{[}{]}{#5}{#7} }
	}
	\NewDocumentCommand{\sqlbrack}{s O{} m}{%
		\IfBooleanTF{#1}
		{ \bracketl*{#3}{[} }
		{ \bracketl[#2]{#3}{[} }%
	}
	\NewDocumentCommand{\sqrbrack}{s O{} m t^G{} t_G{} }{%
		\IfBooleanTF{#1}
		{ \bracketr*{#3}{]}{#5}{#7} }
		{ \bracketr[#2]{#3}{]}{#5}{#7} }%
	}

	\NewDocumentCommand{\cbrack}{s O{} m t^G{} t_G{} }{%
		\IfBooleanTF{#1}
		{ \bracketlr*{#3}{\{}{\}}{#5}{#7} }
		{ \bracketlr[#2]{#3}{\{}{\}}{#5}{#7} }
	}
	\NewDocumentCommand{\clbrack}{s O{} m}{%
		\IfBooleanTF{#1}
		{ \bracketl*{#3}{\{} }
		{ \bracketl[#2]{#3}{\{} }%
	}
	\NewDocumentCommand{\crbrack}{s O{} m t^G{} t_G{} }{%
		\IfBooleanTF{#1}
		{ \bracketr*{#3}{\}}{#5}{#7} }
		{ \bracketr[#2]{#3}{\}}{#5}{#7} }%
	}

	\NewDocumentCommand{\abrack}{s O{} m t^G{} t_G{} }{%
		\IfBooleanTF{#1}
		{ \bracketlr*{#3}{\langle}{\rangle}{#5}{#7} }
		{ \bracketlr[#2]{#3}{\langle}{\rangle}{#5}{#7} }
	}
	\NewDocumentCommand{\albrack}{s O{} m}{%
		\IfBooleanTF{#1}
		{ \bracketl*{#3}{\langle} }
		{ \bracketl[#2]{#3}{\langle} }%
	}
	\NewDocumentCommand{\arbrack}{s O{} m t^G{} t_G{} }{%
		\IfBooleanTF{#1}
		{ \bracketr*{#3}{\rangle}{#5}{#7} }
		{ \bracketr[#2]{#3}{\rangle}{#5}{#7} }%
	}

	\DeclarePairedDelimiter{\abs}{\lvert}{\rvert}

	\NewDocumentCommand{\evalat}{s O{} m m}{%
	  \IfBooleanTF{#1}
	    {\kern-\nulldelimiterspace\left.#3\right|_{#4}}
	    {#3#2|_{#4}}%
	}

	\DeclarePairedDelimiterX{\closedinterval}[2]{[}{]}{#1,\,#2}
	\DeclarePairedDelimiterX{\openinterval}[2]{]}{[}{#1,\,#2}
	\DeclarePairedDelimiterX{\leftopeninterval}[2]{]}{]}{#1,\,#2}
	\DeclarePairedDelimiterX{\rightopeninterval}[2]{[}{[}{#1,\,#2}

	\providecommand\given{}
	\newcommand\SetSymbol[1][]{%
		\nonscript\:#1\vert
		\allowbreak
		\nonscript\:
		\mathopen{}
	}
	\DeclarePairedDelimiterX\Set[1]\{\}{%
		\renewcommand\given{\SetSymbol[\delimsize]}
		#1
	}

	\ExplSyntaxOn
	\NewDocumentCommand{\coord}{sO{}m}
	{
		\IfBooleanTF{#1}
		{\mleft(\coord_print:n {#3}\mright)}
		{\mathopen{#2(}\coord_print:n {#3}\mathclose{#2)}}
	}
	
	\seq_new:N \l_coord_list_seq
	\tl_new:N \l_coord_last_tl
	\cs_new_protected:Npn \coord_print:n #1
	{
		\seq_set_split:Nnn \l_coord_list_seq { , } { #1 } %
		\seq_pop_right:NN \l_coord_list_seq \l_coord_last_tl
		\seq_map_inline:Nn \l_coord_list_seq { ##1 , } %
		\tl_use:N \l_coord_last_tl
	}
	\ExplSyntaxOff

\newcommand{\infrac}[2]{{\begingroup #1 \endgroup /{#2}}}

\newcommand{\Eq}[1]{Eq.~(\ref{#1})}

\newcommand{\Fig}[1]{Fig.~\ref{#1}}

\newcommand{\Sec}[1]{Sec.~\ref{#1}}

\newcommand{\cI}{\mathcal{I}}

\newcommand{\cO}{\mathcal{O}}

\newcommand{\mrm}[1]{\mathrm{#1}}
\newcommand{\eps}{\epsilon}

\newcommand*{\mathmrm}[1]{\relax\ifmmode\mrm{#1}\else{#1}\fi} %

\newcommand{\coshfn}{\cosh\opbrack}
\newcommand{\sinhfn}{\sinh\opbrack}
\newcommand{\tanhfn}{\tanh\opbrack}

\newcommand{\ee}{\mathrm{e}}
\newcommand{\expof}[1]{\ee^{#1}}

\DeclareMathOperator*{\Imaginarypart}{Im}
\newcommand{\Over}[1]{\frac{1}{#1}}

\makeatletter
	\newcommand*{\transpose}{%
		{\mathpalette\@transpose{}}%
	}
	\newcommand*{\@transpose}[2]{%
		\raisebox{\depth}{$\m@th#1\intercal$}%
	}
\makeatother
\newcommand{\dint}[4]{\int_{#1}^{#2}#3\,\mrm{d}#4}

\newcommand{\dd}{\mrm{d}}

\renewcommand*{\vec}[1]{\bm{\mrm{#1}}}

\newcommand{\pitwo}{\frac{\pi}{2}}
\newcommand{\inpitwo}{\infrac{\pi}{2}}

\DeclareMathOperator*{\LambertW}{W}

\newcommand*{\dimension}{D} %
\DeclarePairedDelimiterX{\commutator}[2]{[}{]}{#1,\,#2} %
\DeclarePairedDelimiterX{\acommutator}[2]{\{}{\}}{#1,\,#2} %
\newcommand*{\ret}{\mrm{r}} %
\newcommand*{\adv}{\mrm{a}} %
\newcommand*{\cau}{\mrm{c}} %
\DeclarePairedDelimiterX{\greenf}[2]{\llangle}{\rrangle}{#1;\,#2} %
\DeclarePairedDelimiterXPP{\retgreenf}[2]{}{\llangle}{\rrangle}{_\ret}{#1;\,#2}
\DeclarePairedDelimiterXPP{\advgreenf}[2]{}{\llangle}{\rrangle}{_\adv}{#1;\,#2}
\DeclarePairedDelimiterXPP{\caugreenf}[2]{}{\llangle}{\rrangle}{_\cau}{#1;\,#2}
\DeclarePairedDelimiterXPP{\fougreenf}[2]{}{\llangle}{\rrangle}{_E}{#1;\,#2}
\DeclarePairedDelimiterXPP{\retgreenfeta}[2]{}{\llangle}{\rrangle}{_{\ret,\eta}}{#1;\,#2}
\DeclarePairedDelimiterXPP{\advgreenfeta}[2]{}{\llangle}{\rrangle}{_{\adv,\eta}}{#1;\,#2}
\DeclarePairedDelimiterXPP{\caugreenfeta}[2]{}{\llangle}{\rrangle}{_{\cau,\eta}}{#1;\,#2}

\DeclarePairedDelimiterX{\braket}[2]{\langle}{\rangle}{%
	#1\,\delimsize\vert\,\mathopen{}#2
}
\DeclarePairedDelimiterX{\braopket}[3]{\langle}{\rangle}{%
	#1\,\delimsize\vert\,\mathopen{}#2\,\delimsize\vert\,\mathopen{}#3
}

\newcolumntype{P}[1]{>{\raggedright\arraybackslash}p{#1}} %

\newcommand*{\tstrans}{\Psi} %
\newcommand*{\exactint}{\cI} %
\newcommand*{\discint}[1]{\cI_{#1}} %
\newcommand*{\disctruncint}[2]{Q_{#1}^{#2}} %
\newcommand*{\discerror}[1]{\Delta\discint{#1}} %
\newcommand*{\truncerror}{\epsilon_{\mrm{t}}} %
\newcommand*{\disctruncerror}[2]{\Delta \disctruncint{#1}{#2}} %
\newcommand*{\suminda}{i}
\newcommand*{\sumindb}{j}

\newcommand*{\minfl}{F_{\text{min}}} %
\newcommand*{\minexp}{L} %
\newcommand*{\optimal}{\text{opt}} %
\newcommand*{\maximum}{\text{max}} %
\newcommand*{\evalabsc}{\text{eval}} %
\newcommand{\epsmach}{\eps_{\text{m}}}

\begin{document}
\title{Tanh-sinh quadrature for single and multiple integration using floating-point arithmetic}

\author[1,2]{Joren Vanherck\corref{cor1}}
\ead{joren.vanherck@uantwerpen.be}

\author[1,2,3]{Bart Sor\'ee}
\ead{bart.soree@imec.be}

\author[1,2]{Wim Magnus}
\ead{wim.magnus@uantwerpen.be}

\cortext[cor1]{Corresponding author}
\address[1]{Physics Department, Universiteit Antwerpen, Groenenborgerlaan 171, B-2020 Antwerpen, Belgium}
\address[2]{IMEC, Kapeldreef 75, B-3001 Leuven, Belgium}
\address[3]{ESAT, KU Leuven, Kasteelpark Arenberg 10, B-3001 Leuven, Belgium}

\begin{abstract}
	The problem of estimating single- and multi-dimensional integrals, with or without end-point singularities, is prevalent in all fields of scientific research, and in particular in physics.
	Although tanh-sinh quadrature is known to handle most of these cases excellently,
	its use is not widely spread among physicists.
	Moreover, while most calculations are limited by the use of finite-precision floating-point arithmetic, similar considerations for tanh-sinh quadrature are mostly lacking in literature, where infinite-precision floating-point numbers are often assumed.
	Also little information is available on the application of tanh-sinh quadrature to multiple integration.

	We have investigated the risks and limitations associated with li\-mi\-ted-pre\-ci\-sion floating-point numbers when using tanh-sinh quadrature for both single and multiple integration, while obtaining excellent convergence rates.
	In addition, this paper provides recommendations for a straightforward implementation using limited-precision
	floating-point numbers and for avoiding numerical instabilities.

\end{abstract}


\maketitle

\section{Introduction}
\label{sec::intro}

While many automatic integration routines yield excellent results in general, they often fail in specific circumstances, such as the occurrence of end-point singularities.
In these situations, it is paramount to account for the limitations of a specific routine in order to exploit its full potential.
Tanh-sinh quadrature for numerical integration already exists about half a century,\cite{Takahasi73} but is not well-known in the physics community, though being adopted in some recent publications.\cite{Okayama13,Gaudreau15,Vanherck18}\@
The corresponding scheme can be invoked as an almost general purpose quadrature, which is especially efficient for divergent integrands.
Due to its exponential convergence rate and its good behaviour in general, tanh-sinh quadrature has become rather popular in the field of experimental mathematics.
It even has been coined to be the best scheme for integrands typically encountered in that field, which focuses on high-precision integral calculations.\cite{Bailey05,Bailey11,Slevinsky15}

We describe the quadrature in an accessible manner in \Sec{sec::tanhsinh_quad} while elaborating on its implementation with limited-precision floating-point numbers in \Sec{sec::practicals}.
Specifically, we explain how numerical instabilities can be avoided, while showing how to extend the quadrature to multiple integration, which is indispensable in computational physics.
Examples of integrals calculated with the tanh-sinh rule are given in \Sec{sec::examples}, whereas the conclusions and final remarks are presented in Section~\ref{sec::conclusion}.

\section{Tanh-sinh quadrature}
\label{sec::tanhsinh_quad}

Tanh-sinh quadrature computes integrals of the form%
\footnote{
	Any integral $\dint{a}{b}{f\opbrack{y}}{y}$ with $a$ and $b$ finite can be cast in this form by the linear substitution $2y=\nbrack{b-a}x + \nbrack{b+a}$.
}
\begin{equation} \label{eq:normalizedintegral}
	\cI
	=
		\dint{-1}{1}{ f\opbrack{x} }{x}.
\end{equation}
The quadrature is based on a variable substitution $x=\tstrans\opbrack{t}$, mapping the original, finite domain $x \in \openinterval{-1}{1}$ onto the entire real axis $t \in \openinterval{-\infty}{+\infty}$:
\begin{equation}
	\cI
	=
		\dint{-\infty}{\infty}{
			g\opbrack{ t }
		}{t}
	, \quad
	g\opbrack{ t }
	\coloneqq
		f\opbrack[\big]{ \tstrans\opbrack{t} }
		\tstrans'\opbrack{t}.
\end{equation}
For integration over the entire real axis of exponentially decaying integrands as $\abs{t}\to\infty$, the trapezoidal rule (also known as Sinc quadrature~\cite{Okayama13}) is the most efficient~\cite{Goodwin49,Schwartz69,Sugihara97,Mori05,Trefethen14} among quadratures with equidistant abscissae.
Applying this rule with a step size $h$ between the evaluation points yields the approximate integral
\begin{equation} \label{eq:trapezoidalrule}
	\discint{h}
	=
		\lim\limits_{m\to\infty} \sum_{\suminda=-m}^{+m}
		h \tstrans'\opbrack{t_\suminda}
		f\opbrack[\big]{ \tstrans\opbrack{t_\suminda} },
\end{equation}
the transformed evaluation points being
\begin{equation}
	t_\suminda
	\coloneqq
		\suminda h, \quad i = 0, \pm 1, \pm 2, \ldots
\end{equation}
while introducing a discretization error $\discerror{h}$.
Keeping only $N \coloneqq 2n+1$ function evaluations, we are left with
\begin{equation}
	\disctruncint{h}{n}
	=
		\sum_{\suminda=-n}^{+n}
		h \tstrans'\opbrack{t_\suminda}
		f\opbrack[\big]{ \tstrans\opbrack{t_\suminda} },
\end{equation}
thereby introducing a truncation error $\truncerror$, since the transformed integration domain is now limited to the window $\closedinterval{t_{-n}}{t_{n}}$.
The total error of the resulting approximation as compared to the original integral is limited by $\abs[\big]{ \exactint - \disctruncint{h}{n} }\leqslant \disctruncerror{h}{n} \coloneqq \abs{\discerror{h}} + \abs{\truncerror}$.

\begin{figure}[!tb]
	\centering
	\includegraphics{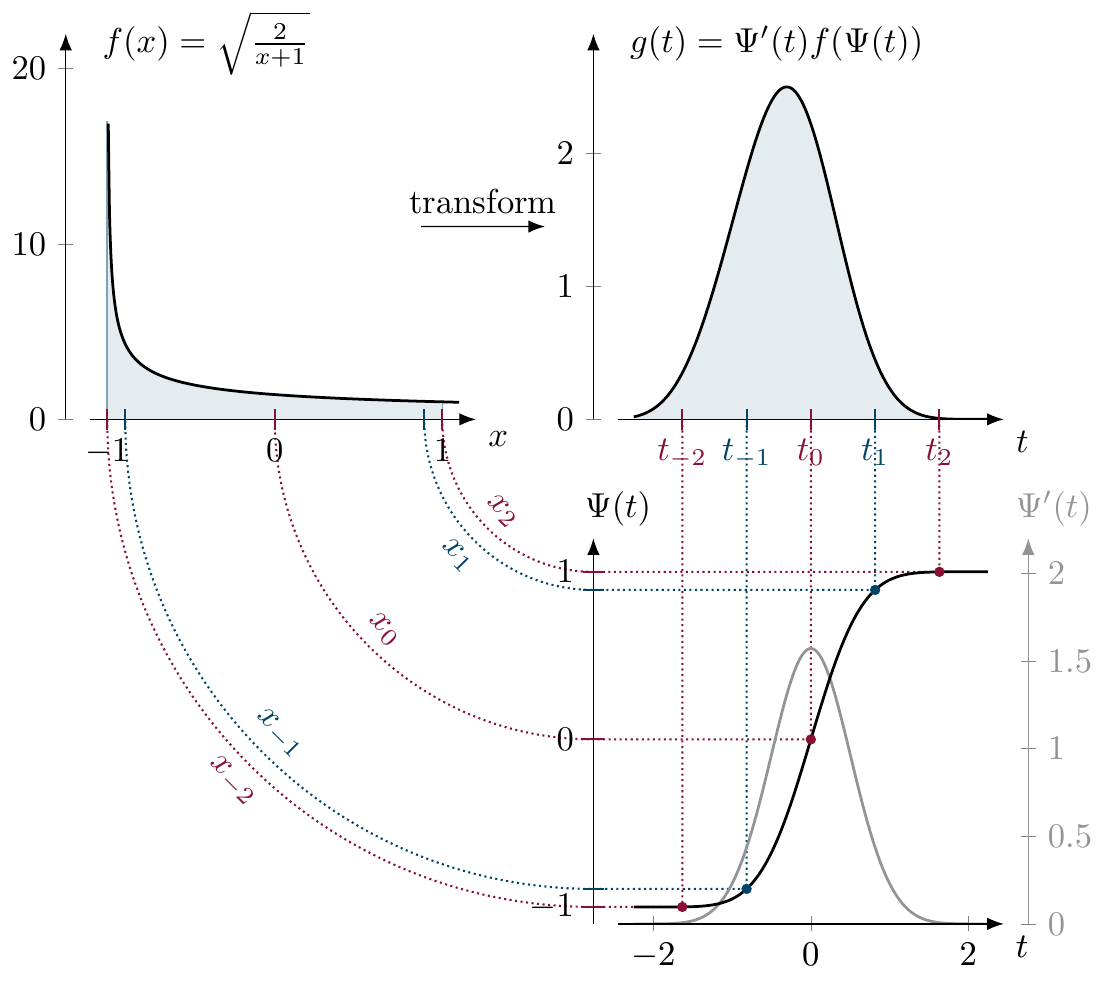}
	\caption{
		Tanh-sinh transformation applied to the integral $\dint{0}{1}{\frac{1}{\sqrt{x}}}{x}$, brought to the standard domain $\openinterval{-1}{1}$ (top left).
		The tanh-sinh transformation transforms the integration domain from $x\in\openinterval{-1}{1}$ to $t\in\openinterval{-\infty}{\infty}$, while the integrands magnitude is strongly suppressed as $t\to\infty$ (e.g.,\ $\evalat{g\opbrack{t}}{t=-2.5} \approx 10^{-7}$ and $\evalat{g\opbrack{t}}{t=-4.5} \approx 10^{-59}$).
		The trapezoidal rule is then applied to this transformed integrand within the window $t\in\closedinterval{-t_n}{t_n}$, for $n=2$ (5 evaluation points) as an illustration (top right).
		The abscissa $x_i$ of the original integration domain to which the equidistant points $t_i$ in the transformed domain correspond, can be recovered as $x_i = \Psi\opbrack{t_i}$ (bottom right and dotted connection lines).
		Notice that most of the abscissae $x_i$ end up very close to the domain boundaries.
		At those points, the effect of the original integrand is strongly suppressed because it is multiplied by $\Psi'\opbrack{t}$ (bottom right, grey).
		This results in the excellent convergence properties of the tanh-sinh quadrature, especially for integrands with singularities near the end-points.
	}
	\label{fig:transform}
\end{figure}

As to the mapping $\Psi(t)$, Takahasi and Mori~\cite{Takahasi73} proposed a \emph{tanh-sinh} transformation
\begin{subequations}
\begin{align}
	\tstrans\opbrack{t}
	&=
		\tanhfn[\big]{\lambda \sinhfn{t}} \label{eq:abscissaeformula}
	\\
	\tstrans'\opbrack{t}
	&=
		\frac{
			\lambda \coshfn{t}
		}{
			\cosh^2\opbrack[\big]{\lambda \sinhfn{t}}
		},\label{eq:weightsformula}
\end{align}
\end{subequations}
where $\lambda=\inpitwo$.
With this choice, which is illustrated in \Fig{fig:transform}, the behaviour of $g\opbrack{ t }$ as $\abs{t}\to\infty$ is optimal in the sense that a faster decay leads to a higher discretization error for a given $h$, while a slower decay causes a larger truncation error because $ \tstrans'\opbrack{t} f\opbrack[\big]{ \tstrans\opbrack{t} }$ remains significant for large values of $\abs{t}$.
An extended error analysis in the complex plane indeed shows that, for a large number of evaluation points, the above transformation often yields an optimal error~\cite{Takahasi73,Sugihara97} as compared to other quadratures based on a variable transformation.
This error is of the order~\cite{Takahasi73,Sugihara97,Mori05,Trefethen14}
\begin{equation} \label{eq:errororder}
	\abs{\disctruncerror{h}{n}}
	=
		\cO\opbrack[\bigg]{
			\exp\opbrack[\Big]{ - \frac{ \pi d N }{ \ln\opbrack{2dN} } }
		}
	,\quad\text{with } N=2n+1
\end{equation}
when $h$ is chosen optimally.
The parameter $d$ has to be chosen such that the analytical continuation of the integrand $g\opbrack{z} = \tstrans'\opbrack{z}f\opbrack[\big]{ \tstrans\opbrack{z} }$ in the complex plane is regular in the strip around the real axis defined by $\abs{\Imaginarypart z} < d$, meaning that the integrand lies within the Hardy space.\cite{Hardy15,Sugihara97}\@
The choice of $d$ and $h$ will be detailed in \Sec{sec::hspacing}, the latter being of particular importance in avoiding numerical instabilities.

The tanh-sinh quadrature scheme has several advantages worth pointing out explicitly.
(i) The error estimate~\eqref{eq:errororder} shows that doubling the number of evaluation points also roughly doubles the number of significant digits, which makes the quadrature especially suited for high precision calculations.\cite{Bailey05,Bailey11}\@
(ii) Even when not all conditions for the error estimate~\eqref{eq:errororder} are strictly fulfilled by the integrand, convergence is still rather fast.
(iii) For integrands with end-point singularities tanh-sinh quadrature converges especially fast as compared to other schemes, thanks to its double-exponential suppression of these divergences (\Fig{fig:transform}).
(iv) The abscissae and weights can be directly extracted from the transformation formulas~\eqref{eq:abscissaeformula} and~\eqref{eq:weightsformula}, in contrast to some other quadrature schemes leaning on iterative processes.
(v) Finally, the transformed integral has equidistant abscissae (trapezoidal rule, \Fig{fig:transform}), which can straightforwardly be reused together with the weights and the function evaluations when doubling the quadrature order is required.

In the present paper we focus on the quadrature scheme as presented above and its generalization to multiple integration, although several extensions are possible.
In general, the transformation function $\Psi\opbrack{t}$ should be chosen such that the transformed integrand $g\opbrack{t}$ behaves like
\begin{equation}
	\abs{g\opbrack{t} }
	\propto
		\exp\opbrack[\Big]{-\pitwo \exp\opbrack{\abs{t}}}
	\quad\text{as}\quad
		t\to \pm \infty
\end{equation}
for the quadrature to be optimal.\cite{Takahasi73}\@
This observation underlies various transformations involving (semi-)infinite domains and integrands with peculiar behavior near the integration boundaries.\cite{Takahasi73}\@
Quadratures based on such transformations are called double exponential (DE) quadratures for obvious reasons.
A second extension is possible when the function $g\opbrack{t}$ has very different decay-rates for $t\to + \infty$ and $t\to -\infty$.
In that case, it can be beneficial to use an unequal number of evaluation points for $t>0$ and $t<0$ when truncating the series of the trapezoid rule~\eqref{eq:trapezoidalrule}.\@
The details, including error estimates with explicit constants, have been worked out by Okayama et al.\ \cite{Okayama13} and more information on asymmetric ranges can be found in Chapter~2 of the book on Sinc methods by Lund and Bowers.\cite{Lund92}\@
Finally, specific modifications for handling oscillatory~\cite{Ooura91,Ooura99} or indefinite integrals are discussed in refs.~\cite{Stenger81,Muhammad03}.

\section{Discussion}
\label{sec::practicals}

The above mentioned favorable properties of tanh-sinh quadrature can only be achieved if numerical instabilities are avoided, especially when dealing with finite-precision floating-point numbers.
Related caveats will be treated and resolved in the following subsection, together with the extension of tanh-sinh quadrature to multi-dimensional integration domains.

\subsection{Avoiding numerical instabilities}
\label{sec::instabilities}

The above presented $N$-point quadrature scheme can be summarized as
\begin{equation}
	\disctruncint{h}{n}
	=
		h\opbrack{n}
		\sum_{\suminda=-n}^{+n}
			w_\suminda
			f\opbrack{ x_\suminda },
\end{equation}
the weights and abscissae respectively being $w_\suminda = \tstrans'\opbrack{t_\suminda}$ and $x_\suminda = \tstrans\opbrack{t_\suminda}$.
It depends crucially on evaluations very close to the end-points of the integration domain, where the density of abscissae $x_\suminda$ is very high (\Fig{fig:transform}).
Correspondingly, one should be careful not to lose significant figures and therefore it is a good practice to store the array~\cite{Takahasi73}
\begin{equation}
	y_\suminda
	=
		1 \pm x_\suminda
	=
		\frac{
			\exp\opbrack[\Big]{\pm \pitwo \sinh{t_\suminda}}
		}{
			\coshfn[\Big]{\pitwo \sinh{t_\suminda}}
		}
	\quad \text{for}\quad
	\suminda \lessgtr 0,
\end{equation}
containing the distances from the respective abscissae $x_\suminda$ to the closest interval bound.

In the context of tanh-sinh quadrature, the most important cause of numerical instability is numerical underflow, occurring when numerical values become smaller then the underflow level (UFL)\@.
This smallest positive normalized floating-point number is $\minfl=2^{\minexp}$, where $\minexp$ is the smallest exponent representable in a given floating-point model.
Specifically, we need to pay attention to the weights, the abscissae and the function evaluations at these abscissae.
Both the smallest weight $w_{n}=\min\Set{w_\suminda}$ and the smallest stored abscissa value $y_n=\min\Set{y_\suminda}$ are determined by window size $t_{n}$, which bounds the transformed integration window $t_{\suminda}\in\closedinterval{-t_{n}}{t_{n}}$.
We will determine a window limit $t_{\maximum}$ and choose $t_{n} \leqslant t_{\maximum}$ to avoid numerical underflow in a given floating-point model.

For the weights, numerical underflow is avoided if $w_n = \tstrans'\opbrack{t_{n}} \geqslant \minfl$, such that the corresponding window limit is
\begin{equation}\label{eq:condition_tw}
	t_{\maximum}^{w}
	=
		\max\Set{t \given \tstrans'\opbrack{t} \geqslant \minfl}
\end{equation}
Similarly, the smallest stored abscissa should exceed the UFL, $y_{n} = 1-x_{n} = 1-\tstrans\opbrack{t_{n}} \geqslant \minfl$, such that we can define
\begin{equation} \label{eq:condition_tx}
	t_{\maximum}^{x}
	=
		\max\Set[\big]{
			t
			\given
				t
				\leqslant
					\tstrans^{-1}\opbrack{1-\minfl}
		}
\end{equation}
as the window limit, where
\begin{equation}
	\tstrans^{-1}\opbrack{1-\minfl}
	=
		\sinh^{-1}\opbrack[\bigg]{
			\frac{
				\ln \opbrack[\big]{ \frac{2}{\minfl}-1 }
			}{
			\pi
			}
		}.
\end{equation}
Both window limits $t_{\maximum}^{w}$ and $t_{\maximum}^{x}$ are \emph{intrinsic} to the used floating-point model as well as to the dimensionality of the integral, as will be discussed in subsection ~\ref{sec::multi_dim}.
Since both conditions always apply simultaneously, it is convenient to introduce the intrinsic window limit $t_{\maximum}^{xw} = \min\Set{t_{\maximum}^{x},t_{\maximum}^{w} }$.
All above mentioned intrinsic quantities are listed in Table~\ref{tbl:maxoptvalues} for a few common floating-point models.
\begin{table}[!htb]
	\centering
	\caption{
		This table shows the smallest representable exponent $\minexp$, the corresponding underflow level $\minfl$, the window limits $t_{\maximum}^{x}$, $t_{\maximum}^{w}$ and $t_{\maximum}^{xw}$, and the corresponding intrinsic maximal optimal order $n_{\maximum}^{xw}$ for various floating-point (fp) models (IEEE 754-2008).
	}
	\begin{tabular}{ r c r@{}l c c c c c }
		\toprule
		fp model
		& $\minexp$	& \multicolumn{2}{c}{$\minfl$} & $\dimension$%
		& 	$t_{\maximum}^{x}$ & 	$t_{\maximum}^{w}$ &	$t_{\maximum}^{xw}$ &	$n_{\maximum}^{xw}$
		\\
		\midrule
		\multirow{2}{*}{single}	& \multirow{2}{*}{-126}	& \multirow{2}{*}{1.} & \multirow{2}{*}{$175\cdot 10^{-38}$}
		& 1,2	& 4.026	& 4.076	& 4.026	& 37
		\\
		&&&& 3	& 4.026	& 3.425	& 3.425	& 18
		\\
		\midrule
		\multirow{2}{*}{double} & \multirow{2}{*}{-1022}	& \multirow{2}{*}{2.} & \multirow{2}{*}{$225\cdot 10^{-308}$}
		& 1,2	& 6.112	& 6.121	& 6.112	& 442
		\\
		&&&& 3	& 6.112	& 5.437	& 5.437	& 201
		\\
		\midrule
		\multirow{2}{*}{extended} & \multirow{2}{*}{-16382}	& \multirow{2}{*}{3.} & \multirow{2}{*}{$362\cdot 10^{-4932}$}
		& 1,2	& 8.885	& 8.886	& 8.885	& 10228
		\\
		&&&& 3	& 8.885	& 8.194	& 8.194	& 4725
		\\
		\bottomrule
	\end{tabular}

	\label{tbl:maxoptvalues}
\end{table}

Besides the calculation of the abscissae and weights, also the integrand evaluations need to be carefully inspected in view of the numerical stability, which can be accomplished by a proper choice of another window limit $t_{\maximum}^{\evalabsc}$.
Examples will be discussed in section \Sec{sec::examples}.

Finally, we take $t_{\maximum} = \min\Set{t_{\maximum}^{xw},t_{\maximum}^{\evalabsc} }$ as the window limit, imposing the maximum on the transformed abscissae values $t_\suminda$ as to avoid numerical instabilities.

\subsection{Multiple integration}
\label{sec::multi_dim}

Whereas higher dimensional integrals ($\dimension > 1$) are commonplace in physics,~\cite{Bailey11,Slevinsky15} we restrict ourselves to integration domains that are Cartesian products of lower dimensional domains, or that can be transformed to such a region by variable substitutions or coordinate transformations.
In these cases, it is possible to apply the one-dimensional quadrature rules repeatedly in each dimension, but, once again, caution is in order to avoid numerical underflow.
Slightly more complex multiple integrals, having boundaries that are integration variables in an exterior loop, are accessible through the use of indefinite integration.\cite{Muhammad05}

As an example in two dimensions, consider the integration of a function $f\opbrack{x,y}$ over the domain $-1 \leqslant x, y \leqslant 1$, which is typically approximated as
\begin{equation} \label{eq:multi_dimensional}
	\dint{-1}{1}{
		\dint{-1}{1}{
			f\opbrack{x,y}
		}{y}
	}{x}
	\approx
		\disctruncint{h_1,h_2}{n_1,n_2}
	=
		h_1 h_2
		\sum_{\suminda=-n_1}^{+n_1}
		\sum_{\sumindb=-n_2}^{+n_2}
			w_\suminda
			w_\sumindb
			f\opbrack{ x_\suminda, x_\sumindb },
\end{equation}
where $n_1\neq n_2$ in general, although $n_1=n_2$ is often a very convenient choice.
While it is clear from Eq.~\eqref{eq:multi_dimensional} how to extend the quadrature rule to even higher dimensions, it is recommended to calculate the weights $w_\suminda=\Psi'\opbrack{t_{\suminda}}$ in each dimension separately, and multiply them with the function evaluation on the fly, as underflow problems may quickly arise otherwise.
It is reasonable to limit the weight-specific window limit $t_{\maximum}^{w}$, for which we take the same value in each direction, even more in higher dimensions (multiplying many very small weight quickly renders all the weights effectively zero).
A good value for $t_{\maximum}^{w}$ was found to be such that
\begin{equation}
	t_{\maximum}^{w} = \max \Set[\big]{t \given \nbrack{\tstrans'\opbrack{t}}^{\mathcal{D}} \geqslant \minfl }\text{, where }
		\mathcal{D}
		=
		\max\Set[\big]{1,\dimension-1},
\end{equation}
since then $w_n^\mathcal{D}\geqslant \minfl$.
This ensures that each weight separately exceeds the UFL\@.
For higher dimensions, the product of all but one weight should also be higher then the UFL in order that the corresponding term contribute to the sum.
In the latter case, one weight can be disregarded, because it can be compensated by a large value from the function evaluation.
Examples for one, two and three dimensions are given in Table~\ref{tbl:maxoptvalues}.

\subsection{Abscissae spacing}
\label{sec::hspacing}

For most quadrature schemes, the choice of the quadrature order $n$, together with a window size that is fixed by the integration bounds, determines the abscissae spacing $h$.
In the tanh-sinh quadrature schemes, the window size $\closedinterval{-t_{n}}{t_{n}}$ is not fixed a priori, with bounds $t_{n}=n \cdot h$ that depend themselves on the the abscissae spacing.
A smart choice for $h\opbrack{n}$ is indispensable.
We will discuss two alternatives: (i) the optimal $h_{\optimal}\opbrack{n}$ and (ii) the \emph{maximal} $h_{\maximum}\opbrack{n}$.

The optimal spacing
\begin{equation}
	\label{eq::h_opt_d}
	h\opbrack{n}
	=
	h_{\optimal}\opbrack{n}
	\coloneqq
	\frac{2}{N} \LambertW\opbrack{2 d N}
	,\quad\text{with } N=2n+1
\end{equation}
leads to the optimal error bounds~\eqref{eq:errororder} by making the discretization error $\discerror{h}$ and the truncation error $\truncerror$ contribute equally.\cite{Takahasi73,Sugihara97,Mori05,Trefethen14}\@
The Lambert W-function $\LambertW\opbrack{z}$ is implicitly defined as the solution of $z=w\expof{w}$ for $w$.
One should note here that, in most of literature, the optimal width is given by $h = \infrac{2}{N} \ln\opbrack{2 d N}$, which is only correct in the limit of large $N$.\cite{Trefethen14}\@
The optimal abscissae spacing~\eqref{eq::h_opt_d} depends on the strip width $d$ of regularity of the transformed integrand around the real axis.
Strictly speaking, this strip width needs to be specifically determined for each integral to achieve optimal convergence.
However, when the integrand is the result of a complex numerical routine, such a determination of $d$ is not always possible, let alone desirable when a large range of integrals needs to be calculated.
In practice, $d=\infrac{\pi}{2}$ often leads to reasonable convergence rates.
The optimal transformed abscissae spacing then becomes
\begin{equation} \label{eq::h_opt}
	h_{\optimal}\opbrack{n}
	=
		\frac{2}{N} \LambertW\opbrack{\pi N}
	,\quad\text{with } N=2n+1.
\end{equation}

While $h_{\optimal}\opbrack{n}$ is theoretically optimal in an infinite-precision context, its outstanding results can not always be achieved due to the limitations of floating point types.
The main problem is the window size $t_{n}= n h_{\optimal}\opbrack{n}$ that can quickly reach beyond the window limit $t_{\maximum}$, sometimes before full convergence is reached.
In order to avoid numerical instabilities, the optimal spacing should never be used for orders higher than $n_{\maximum}=\max\Set{n \given n h_{\optimal}\opbrack{n} \leqslant t_{\maximum}}$ at which this happens.
This maximal order can depend on numerical difficulties in the integrand evaluation through $t_{\maximum}$ and must therefore be determined on a per-case basis.
The maximal order itself is always limited by the intrinsic maximal order $n_{\maximum}^{xw}$ which is fully determined by the floating-point type and dimension (Table~\ref{tbl:maxoptvalues}).

As an alternative to the optimal spacing, we propose to use the maximal spacing
\begin{equation} \label{eq::h_max}
	h\opbrack{n}
	=
		h_{\maximum}\opbrack{n}
	\coloneqq
		\frac{t_{\maximum}}{n}
\end{equation}
instead.
Here, the window size $t_{n}$ is fixed, independent of the order $n$, to $t_{\maximum}$ so as to avoid numerical underflow.
Besides the obvious simplicity of the $h_{\maximum}\opbrack{n}$, also the reusability of function evaluations obtained from lower order estimates, is an advantage when it comes to doubling the order index $n$.

\section{Worked examples}
\label{sec::examples}

In this section, we demonstrate the use of the tanh-sinh quadrature as described above, specifically comparing results using optimal and maximal spacing for different floating-point types.
For each example integral, possible pitfalls and noteworthy features are discussed, while showing the relative error as a function of the quadrature order.
This relative error is intrinsically limited from below by the machine precision $\epsmach$, which is indicated by dashed lines for each floating-point models.
The results for optimal spacing can only be obtained up to the maximum quadrature order $n_{\maximum{}}$, which is shown as a larger data point if it falls within the shown domain.
Results obtained using double-precision Gauss-Legendre quadrature with the same number of evaluation points are always shown for comparison.
We focus on integrands with a singularity near the origin of the original integration domain, as may be relevant to solve numerous problems in physics.
For benchmarking purposes, it is not the intent to evaluate the specific integral, but rather to demonstrate how the quadrature behaves for different types of limiting behavior.

\subsection{Integrands with a non-integrable singularity}

\begin{figure}[!tb]
	\centering
	\includegraphics{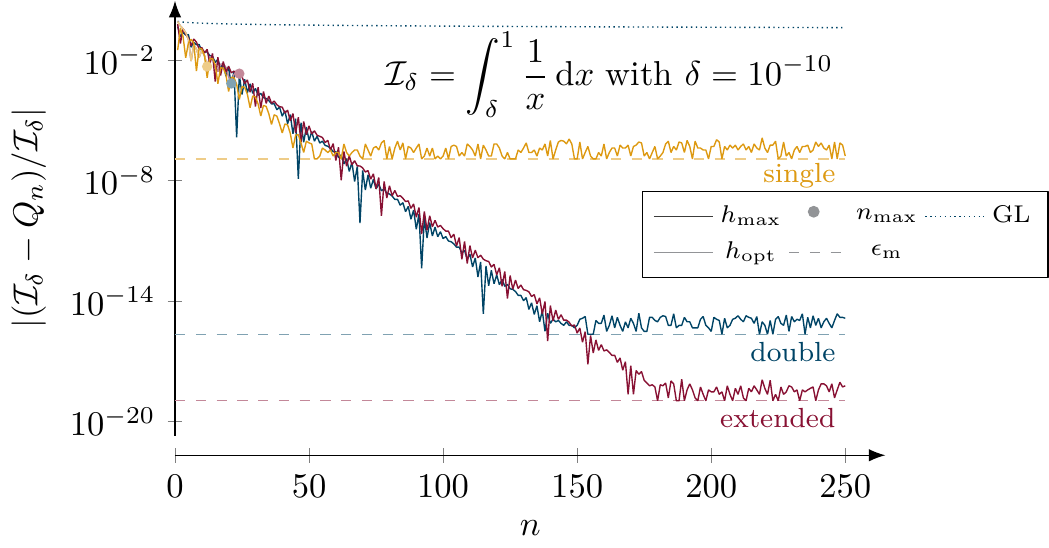}
	\caption{
		The relative error decreases exponentially with quadrature order $n$.
		Both optimal and maximal spacing yield similar convergence rates.
		However, whereas the maximal spacing reaches full precision for all floating-point types, the optimal spacing can barely converge because $n_{\maximum}$ is very small.
		The fluctuations as a function of $n$ decrease upon increasing $a$.
		No convergence can be distinguished when using Gauss-Legendre quadrature.
	}
	\label{fig:epserror}
\end{figure}

Sometimes, one needs to integrate a function with a non-integrable singularity located close to one of the integration limits.
A straightforward example of this type is the integral
\begin{equation} \label{eq:physics:overeps}
	\cI_{\delta}
	=
		\dint{\delta}{1}{ \Over{x} }{x}
	=
		-\ln \delta,
\end{equation}
where $\delta$ is an arbitrarily small positive number.
In spite of $\ln \delta$ quickly diverging when $\delta$ tends to zero, tanh-sinh quadrature quickly yields remarkably high precision integrals when the integrand is sampled with many points near $x=\delta$ (\Fig{fig:epserror}).
This can, however, also be a source of errors: if the sampling points are closer to $x=\delta$ than machine precision permits, the function evaluations will be inaccurate.
To prevent this, $t_{\maximum}^{\evalabsc}$ should be chosen such that the abscissae $x=x_{-n}$ remains distinguishable from $x=\delta$.
For the errors shown in \Fig{fig:epserror}, we chose $t_{\maximum}^{\evalabsc}$ such that $\infrac{\nbrack{x_{-n} - \delta}}{\delta} > a\cdot\epsmach$, with $a=100$.
Increasing the value of $a$ diminishes the fluctuations of the resulting error, but increases the achieved minimal relative error.

Notice in \Fig{fig:epserror} that $n_{\maximum}$ is very small and the corresponding error using optimal spacing is still far above the achievable error.
The use of the proposed maximal spacing does allow for a steady and fast decrease of the relative error up to full precision.
Gauss-Legendre quadrature fails completely in this test.
It lacks the fine-grained sampling around $x=\delta$.

\subsection{Multiple integration with singular integrands}

\begin{figure}[!tb]
	\begin{adjustwidth}{-3cm}{3cm}
	\centering
	\includegraphics{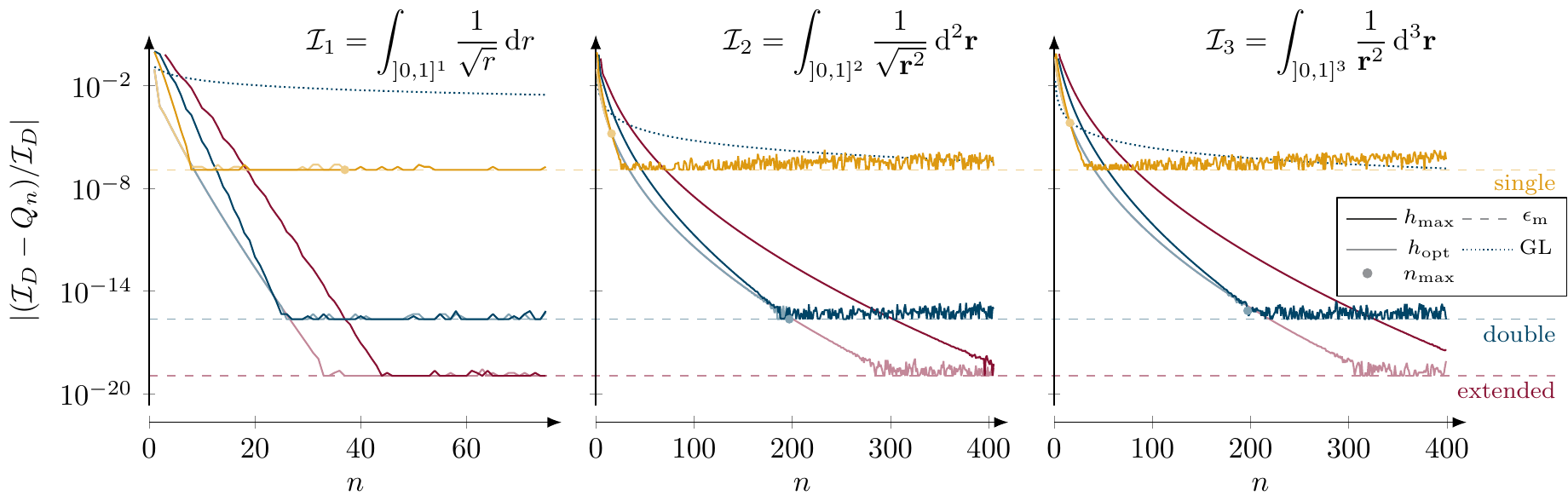}
	\end{adjustwidth}
	\caption{
		The relative error as compared to the exact result for one-dimensional quadrature order $n$ for functions with an integrable singularity at the origin in one, two and three dimensions.
		The total number of sampling points in $\dimension$ dimensions is $(2n+1)^{\dimension}$.
		The quadrature with optimal spacing always converges slightly faster than the maximal spacing variant.
		However, it does not reach full precision for the single-precision floating-point type in the two- and three-dimensional examples, whereas the quadrature with $h_{\maximum}\opbrack{n}$ does.
		While still being fast, the convergence is more slow in higher dimensions.
		Gauss-Legendre quadrature converges, but only very slowly.
	}
	\label{fig:divergence_fdim}
\end{figure}

Next, we consider $\dimension$-dimensional integrals
\begin{equation}
	\cI_{\dimension}
	=
		\int\limits_{\mathclap{\leftopeninterval{0}{1}^\dimension}}^{} \dd^{\dimension}\vec{r}\,f_{\dimension}\opbrack{\vec{r}},
\end{equation}
where the integrands are chosen to be
\begin{align*}
	f_1\opbrack{x}
	&=
		\frac{1}{\sqrt{x}}
	\\
	f_2\opbrack{x,y}
	&=
		\frac{1}{\sqrt{x^2+y^2}}
	\\
	f_3\opbrack{x,y,z}
	&=
		\frac{1}{x^2+y^2+z^2},
\end{align*}
all containing a integrable singularity at $\vec{r}=0$.
Integrals of the type $\cI_{\dimension}$ are typically encountered in solid state physics, where integration over the first Brillouin zone is very common.
For example, in a recent study of anisotropic quantum Heisenberg ferromagnets~\cite{Vanherck18}, one of the integrands that determines the Curie temperature falls of in the same way as $f_3(\vec{r})$.
The results\footnote{%
	Elementary integration of the function class $f_{\dimension}\opbrack{\vec{r}}$ leads to $\cI_1 = 2$, $\cI_2 = 2 \ln(1 + \sqrt{2})$ and $\cI_3 = 3 \nbrack{ \operatorname{Ti_2}\opbrack{3-2\sqrt{2}} - C} + \frac{3\pi}{4} \tanh^{-1}\opbrack{\frac{2 \sqrt{2}}{3}}$, where $\operatorname{Ti_2}\opbrack{x}$ and $C$ are respectively the inverse tangent integral and the Catalan constant.
} for each of these integrals are shown in \Fig{fig:divergence_fdim}.
The convergence rate is especially large for the one-dimensional integral, but also for higher dimensions the convergence is fast and steady up to maximum precision.
While no special numerical concerns apply to $f_1$, numerical underflow in the calculation of the denominator of $f_2$ and $f_3$ needs to be avoided by choosing $t_{\evalabsc}$ such that $x_{-n} \geqslant \sqrt{\minfl}$.

In the latter two cases, the maximal order $n_{\maximum{}}$ was reached in the single precision optimal scheme before convergence up to full precision was achieved.
In the double-precision equivalent, full precision is reached at an order just below $n_{\maximum{}}$.
All convergence rates seem to become lower when the dimensionality of the integral is increased.
This might be due to the fact that the choice of $d=\inpitwo{}$ in \Eq{eq::h_opt_d} is not optimal for the outer integrals.
The scheme using maximal spacing always reaches full precision.
Even though its convergence rate is somewhat slower than the optimal (especially for the extended precision calculation), it still has a decent performance.
In every case, Gauss-Legendre quadrature converges, be it at a much lower rate than the tanh-sinh quadratures.

\section{Conclusion}
\label{sec::conclusion}

We have reviewed tanh-sinh quadrature, which is very efficient for integrands with end-point singularities, that are often encountered in physics.
We demonstrated that even limited-precision floating-point arithmetic facilitates fast convergence at machine precision levels.
Extending the scheme to higher dimensions turns out to be straightforward, although the convergence rate is typically lower as a function of the quadrature order.

We showed that care has to be taken to avoid numerical underflow in the abscissae and the weights, as well as to avoid other numerical instabilities related to integrand evaluations.
These problems can typically be resolved straightaway by choosing a suitable window limit $t_{\maximum}$.
For the optimal spacing $h_{\optimal{}}\opbrack{n}$, this implies a maximum usable quadrature order $n_{\maximum}$, sometimes hampering convergence.
The proposed maximal spacing rule $h_{\maximum}\opbrack{n}$ alleviates this limitation.
Despite having slightly slower convergence rates, this rule achieves full machine precision accuracy more consistently.
Moreover, its implementation is even more straightforward and it allows for the reuse of abscissae, weights and function evaluations when it comes to recompute a particular integral for increasing quadrature orders, as is often convenient to monitor convergence.

\end{document}